\documentclass[10pt]{article}

\usepackage{graphicx}
\usepackage{amssymb,amsmath}
\usepackage{latexsym}
\usepackage{algorithm}
\usepackage{algpseudocode}

\newcommand{\beq}{\begin{equation}}
\newcommand{\eeq}{\end{equation}}
\newcommand{\beqa}{\begin{eqnarray}}
\newcommand{\eeqa}{\end{eqnarray}}
\newcommand{\beqas}{\begin{eqnarray*}}
\newcommand{\eeqas}{\end{eqnarray*}}
\newcommand{\bit}{\begin{itemize}}
\newcommand{\eit}{\end{itemize}}
\newcommand{\benum}{\begin{enumerate}}
\newcommand{\eenum}{\end{enumerate}}

\newcommand{\mmp}[1]{{\em MMP: #1}}
\newcommand{\comment}[1]{}

\newcommand{\xent}{\ensuremath{\operatorname{XE}}}
\newcommand{\gmms}{\ensuremath{GMM}}
\newcommand{\mm}{{\rm MM}}

\newcommand{\gmmsup}{\ensuremath{^{GMM}}}  

\newcommand{\bigOO}{{\cal O}}

\newcommand{\bbone}{\mathbf{1}}

\newcommand{\ssss}{{\mathbb S}}  
\newcommand{\expe}{{\mathbb E}}  
\newcommand{\pimin}{\pi^{-1}}
\newcommand{\sigmin}{\sigma^{-1}}
\newcommand{\pizero}{{\sigma}}

\newcommand{\idperm}{{\rm id}}

\newcommand{\uppper}{^{\rm up}}

\newcommand{\listab}{{\tt list}}
\newcommand{\pijtab}{{\tt Pij}}

\newcommand{\rank}{\operatorname{rank}}
\newcommand{\geom}{\operatorname{geom}}

\newcommand{\E}{{\mathcal E}}   
\newcommand{\tilE}{\tilde{\mathcal E}}   
\newcommand{\sbar}{\bar{\mathbf s}}   
\newcommand{\qbar}{\bar{\mathbf Q}}   
\newcommand{\qbarl}{\qbar\uppper}

\newcommand{\myprec}{\prec}
\newcommand{\mysucc}{\succ}
\newcommand{\thetavec}{\vec{\theta}}
\newcommand{\alphavec}{\vec{\alpha}}

\newtheorem{prop}{Proposition}
\newtheorem{corollary}[prop]{Corollary}
\newtheorem{example}{Example}

\title{The Entropy and Crossentropy of Generalized Mallows Models}
\author{Marina Meil\u{a}\\
  Department of Statistics\\
  University of Washington\\
mmp@stat.washington.edu}

\begin{document}

\maketitle

\begin{abstract}
The Generalized Mallows Model (GMM) is a well known family of models for ranking data. A GMM is a distribution over $\ssss_n$, the set of permutations of $n$ objects, characterized by a {\em location parameter} $\pizero\in\ssss_n$, known as {\em central permutation} and a set of dispersion parameters $\theta_{1:n-1}\in (0,1]$. The GMM shares many properties, such as having sufficient statistics, with exponential models, thus it can be seen as an exponential family with a discrete parameter $\pizero$. This paper shows that computing entropy, cross-entropy and Kullback-Leibler divergence in the the class of GMM is tractable, paving the way for a better understaning of this exponential family.
\end{abstract}

\comment{
\cite{fligner:86}
\cite{MPhadnisPattersonBilmes:uai07}
\cite{Stanley:97}
\cite{mallows:57}
\cite{fligner:86}
\cite{Pearl:84}
\cite{MandhaniM:heuristics-aistats09}
\cite{fligner:88}
\cite{cohen:99}
\cite{fukunaga:75}
\cite{cheng:95}
\cite{Critchlow:85}
\cite{Degroot:75}
\cite{jester}
\cite{lebanon:08}
\cite{gormley:06}
\cite{fligner:90}
\cite{quadrianto:09}
\cite{lebanon:03}
\cite{Critchlow:85}
\cite{murphy:03}
\cite{murphy:05}
\cite{gormley:06}
\cite{plackett:75}
\cite{Luce:59}
\cite{murphy:05}
\cite{busse:07}
\cite{MBao:uai08-infinite}
\cite{VanderVaart:98}
}

\section{Introduction}
\label{sec:intro}

Mallows-type models have recently received much interest in the
literature on modeling rankings, and a commensurately large number of
applications. Models in this class have been studied, among others, by
[1],[2], and recently by [3].  Many of the
desirable properties of Mallows models derive from the fact that these
are exponential families where the sufficient statistics are inversion
counts. This means the probability of a permutation $\pi$ is penalized
exponentially in the number of inversions w.r.t a central permutation
$\pizero$ [2, 3]. Moreover, the normalization constant $Z$ can be computed in
closed form and depends only on the dispersion parameters and not on
the location parameter $\pizero$.

Due to these attractive properties, Maximum
Likelihood (ML) estimation of the GMM family were extensively studied 
\cite{fligner:86,fligner:88,MPhadnisPattersonBilmes:uai07,MBao:jrss-infinite07,MArora:13} and efficient estimation algorithms exist 
\cite{MPhadnisPattersonBilmes:uai07,MandhaniM:heuristics-aistats09,WagnerM:RIM22}.  \footnote{It is knows that 
estimating $\pizero$ is NP-hard \cite{bartholdi:89} in general, but in \cite{MPhadnisPattersonBilmes:uai07} it was shown that when consensus exists, the problem can become tractable. For example, in \cite{WagnerM:RIM22} a stochastic
algorithm was introduced that searches the space of permutations
efficiently, and whose every iteration is tractable.}.

Here we turn to the information theoretic properties of the
\gmms~family, such as entropy, cross-entropy, and Kullbach-Leibler
(KL) divergence. This is a necessary step towards elucidating the
properties of this family, as an exponential family with a discrete
parameter $\pizero$.

The paper is structured as follows. I the next section, Section
\ref{sec:bg}, we introduce the \gmms~family and some of its
remarkable properties, which stem from the existence of sufficient
statistics, in the form of the inversion matrix $Q$. Section
\ref{sec:crossent} is the main statistical contribution, deriving
compact expressions for the crossentropy, entropy and KL-divergence,
that align with the known expressions for exponential families, even
though this family has a permutation $\pizero$ as location parameter.
Thus, the crossentropy depends on cross-expectations of sufficient
statistics, averages over $n!$ permutations, which at first sight may
appear intractable.  Section \ref{sec:qbar} contains the more
surprising contribution of the paper, a polynomial time recursion to
evaluate the required expectations.
Examples and experimental illustrations are presented in section \ref{sec:exe} and Section \ref{sec:discussion} concludes the paper.

\section{Background}
\label{sec:bg}
\subsection{Statistical Models for Permutations}
\label{sec:bg-gmm}

Consider a set of $n$ items, which we will denote $\E=\{e_1,e_2,...,e_n\}$, which might represent products in a store, candidates in an election, or statements in a survey. A ranker (typically a person, natural process, or automated system) produces a ranking of the elements in the set, according to their preference (e.g. the best candidate first, the statement they most agree with first). The set of all permutations of $\E$ is denoted $\ssss_\E$, or simply $\ssss_n$, by identifying $e_i$ with $i=1,\ldots n$. By contrast, responses on a Likert scale represent {\em ratings}, another modality of expressing preferences.

A {\em ranking model} is a distribution over $\ssss_n$. This is a
discrete distribution over the $n!$ rankings, far too large a sample
space for most data sets to assign the probability of each possible
ranking by enumeration. This large sample space necessitates models
that are either simply parameterized or highly structured, such as the
Mallows Model \cite{mallows:57}, Generalized Mallows Model
\cite{fligner:86}, and Bradley-Terry Model \cite{BradleyTerry},
Recursive Inversion Model \cite{WagnerM:RIM22}, Plackett-Luce Model
\cite{Luce:59,plackett:75}, and Riffle Independence Model
\cite{huang2012}. The models we develop rely on these same principles.

In this paper we study the Generalized Mallows model, to be described below, where the probability of a permutation $\pi\in\ssss_n$ depends on the inversions between $\pi$ and a modal reference permutation $\pizero$. 

\subsection{The Inversion Matrix $Q$}
We start with establishing basic notation. We represent a permutation
$\pi$ of $\E=\{e_1,e_2,...,e_n\}$ as an ordered list
$\pi=[\pi_1,\,\ldots \pi_n]$, where we denote by $\pi_i$ the item at
rank $i$, and by $\pi^{-1}(e)$ the rank of item $e\in\E$ in $\pi$. For
two items $e,e'\in\E$, when $\pi^{-1}(e)<\pi^{-1}(e')$, meaning that
item $e$ is ranked before item $e'$ in $\pi$, we say that $e$ {\em
  precedes} $e'$ in $\pi$, and write it alternatively as $e\myprec_\pi
e'$. Let $\idperm$ denote the {\em identity permutation} of $\E$; when $\E=\{1,2,...,n\}$ we have $\idperm=[\,1\,\ldots \,n\,]$.

An {\em inversion} between two permutations $\pi$ and $\pizero$ is a
pair of items $e,e'\in\E$, $e\neq e'$, so that $e\myprec_\pi e'$ and
$e \mysucc_{\pizero}e'$. This definition is symmetric
w.r.t. $\pi,\pizero$. In this paper, and more generally in
the ranking models context, $\pizero$ typically represents a {\em
	reference permutation}, such as the identity. The {\em inversion distance} between two permutations $\pi,\pizero$ is the number of inversions between them. 
\beq\label{eq:mallow-kendall}
d(\pi,\pizero)=\sum_{e\in\E}\sum_{e'\in\E}1_{[e \mysucc_{\pizero}e']}\cdot 1_{[e\myprec_\pi e']}.
\eeq 
The inversion distance of $\pi$ w.r.t. $\idperm$ (or another $\pizero$) is also called the {\em number of inversions} of $\pi$ (w.r.t. $\pizero$), or the {\em Kendall tau} distance \cite{fligner:88}. 

Any permutation $\pi$ can be uniquely represented by the  {\em inversion matrix} $Q(\pi)$, which contains information on inversions between items in $\pi$ w.r.t. the identity $\idperm$. 
\beq \label{eq:qpi}
Q(\pi)\;=\;[Q_{ee'}]_{e,e'\in\E},
\quad\,Q_{ee'}=1\;\text{if}\,e\myprec_\pi e'\text{ and } 0
\text{ otherwise, for all }e,e'\in \E;
\eeq
the rows and columns of $Q$ are ordered by $\idperm$.

It is easy to see that the lower triangle of $Q$ is sufficient to reconstruct $Q$, since  $Q_{ee'}=1-Q_{e'e}$. Moreover, $\pi$ can be easily reconstructed from $Q(\pi)$ by summing the columns of $Q$. Indeed, for any $e\in\E$, $\pimin(e)=1+\sum_{e'\in \E}Q_{e'e}$.

The inversion distance $d(\pi,\idperm)$ can be obtained from $Q(\pi)$
as the sum of the elements in its lower triangle. More generally, the
inversion distance between two permutations $\pi,\sigma$ is the sum of the lower triangle of $Q(\sigma)$ when the rows and columns are sorted by $\pi$, or equivalently
\beq
\label{eq:d-pisigma}
d(\pi,\sigma)\;=\;d(\sigma,\pi)\;=\;
d(\pi\circ\sigmin, \idperm) \;=\;
\sum_{e'\in\E}\sum_{e\myprec_\pizero
  e'}Q_{e'e}.  \eeq

\begin{example}\label{ex:mallow}
Consider the permutations $\pi=[c, a, d, b, f, e]$, and $\idperm=[a,b,c,d,e,f]$. We can calculate the number of inversions of $\pi$ easily from the inversion matrix $Q$, shown below.
\beq\label{eq:mallow-q}
Q(\pi;\pizero)\;=\;
\begin{array}{c|c|c|c|c|c|c|c|}
	& a & b & c & d & e & f \\
	\hline
	a & - & 1 & 0 & 1 & 1 & 1\\
	b & 0 & - & 0 & 0 & 1 & 1\\
	c & 1 & 1 & - & 1 & 1 & 1\\
	d & 0 & 1 & 0 & - & 1 & 1\\
	e & 0 & 0 & 0 & 0 & - & 0\\
	f & 0 & 0 & 0 & 0 & 1 & -\\
	\hline
\end{array}
\eeq
In the above, $Q_{a,b}=1$, implying that $a$ comes before $b$, while
$Q_{a,c}=0$ indicating that $c$ precedes $a$. The pairs of items
$(a,c)$, $(b,c)$, $(b,d)$, and $(e,f)$ are out of order in $\pi$, matching the four values of 1 found in the
lower triangle of the inversion matrix. Thus, $d(\pi,\idperm)=4$.

The columns sums of $Q(\pi)$, are 1, 3, 0, 2, 5, 4; therefore, the ranks are $\pimin(a)=1+1=2$, $\pimin(b)=3+1=4$, $\pimin(c)=0+1=1$, etc.
\end{example}

For more details on inversion matrices and their properties relevant
to ranking models the reader is invited to
consult previous works \cite{stanley:Enum,MPhadnisPattersonBilmes:uai07}. As will be shown,
inversion matrices serve as sufficient statistics for the Mallows family of models.

\subsection{The Mallows Model (\mm)}
Exponential ranking models seek to define a probability for each $\pi$, w.r.t. some central (typically modal) ranking denoted as $\pizero$. The simplest of these models is the Mallows model.

The Mallows model \cite{mallows:57} is defined by a 
modal ranking $\pizero$, and a single parameter $\theta\in[0,1]$ which
represents the probability of an inversion w.r.t. the modal
ranking. Since the total number of inversions is the Kendall's tau
distance, the probability of a permutation $\pi$  follows an
exponential distribution, i.e.
\beq\label{eq:mallow-like}
P^\mm_{\pizero,\theta}(\pi)=\frac{1}{Z(\theta)}\theta^{d(\pi,\pizero)},
\eeq
where $Z(\theta)$, is the normalization constant and can be easily derived as
\beq\label{eq:mallow-z}
Z(\theta)=\prod_{k=0}^{n-1}\sum_{i=0}^k\theta^k=(1)(1+\theta)(1+\theta+\theta^2)...(1+\theta+...+\theta^{n-1})
\eeq
which can be computed in closed form \cite{mallows:57}. This is one of the most remarkable properties of Mallows and Generalized Mallows models. 

It can be observed from the probability in
equation~\eqref{eq:mallow-like}, under this parameterization, that $\theta=1$ represents the uniform distribution over $\ssss_n$, values of $\theta$ near 1 correspond to models with large dispersion, while values of $\theta$ near 0 represent distributions concentrated around $\pizero$.

\begin{example} 
For a Mallows model parameterized by $\pizero=[a,b,c,d,e,f]$ and $\theta$,
the probability of the ranking $\pi$ from Example~\ref{eq:mallow-q} is
\beq\label{eq:mallow-example}
P^{\mm}_{\pizero,\theta}(\pi)=\frac{1}{Z(\theta)}\theta^{4}
\eeq
and $P(\pizero|\pizero,\theta)=\frac{1}{Z(\theta)}$ with $Z(\theta)=(1+\theta+\theta^2+\theta^3+\theta^4+\theta^5)\ldots (1+\theta+\theta^2)(1+\theta)$
\end{example}

\subsection{Generalized Mallows Models (GMM)}
The Mallows model penalizes each inversion equally by a factor $\theta$.
This model can be expanded in an elegant way to improve its flexibility,
by replacing the single parameter $\theta$ with
a vector $\thetavec\in[0,1]^{n-1}$. while still utilizing a modal
permutation $\pizero$. The {\em Generalized Mallows Models} \cite{fligner:88}
is parameterized by $\pizero$ and $\thetavec=(\theta_1,\ldots \theta_{n-1})$. We decompose $d(\pi,\pizero)=\sum_{i=1}^{n-1}s_r(\pi|\pizero)$, where each $s_r$ counts inversions w.r.t. rank $r$ of $\pi$.

This decomposition can be best understood as a stagewise construction of
$\pi$ from $\pizero$. One starts from an empty $\pi$; at stage
$r=1,2,\ldots n$, one samples
$s_r\sim \geom(\theta_r; n-r)$, and selects the item on rank $1+s_r$ of
$\pizero$  (alternatively this can be viewed as skipping over $s_r$ items in $\pizero$). This item becomes  $\pi_r$ and is deleted from $\sigma$. 

Consequently, the probability of a permutation $\pi$ is
\beq\label{eq:gmms-like}
P^{\gmms}(\pi|\pizero,\thetavec)=\prod^{n-1}_{r=1}\frac{1}{Z_{n-r+1}(\theta_r)}\theta_r^{s_r(\pi|\pizero)}
\eeq
where again $Z_{n-r+1}(\theta_r)$ represents the normalization constant of the respective geometric distribution. 

In has been demonstrated \cite{MBao:jrss-infinite07} that the $s_r(\pi|\pizero)$ terms can be calculated, like the Kendall's tau distance, from the inversion matrix $Q(\pi,\pizero)$ by
\beq
s_r(\pi|\pizero)\;=\;\sum_{e'\myprec_{\pizero}e} Q_{ee'}(\pi|\pizero)
\quad\text{where }e=\pi(r).
\eeq
Hence $s_r$ are row sums, in which only entries in the lower triangle of $Q$ are counted.
\begin{example}
Let $\pi=[c, a, d, b, f, e]$ as before, and a \gmms~parameterized by
$\pizero=[a,b,c,d,e,f]$ and $\thetavec=(\theta_1,\theta_2,\theta_3,
\theta_4,\theta_5)$. Refer back to equation~\eqref{eq:mallow-q} for the matrix $Q(\pi;\pizero)$. 

From this matrix we can calculate $s_{1:5}(\pi|\pizero)$ as follows: $\pi_1=c$, hence $s_1=\sum_{e'\myprec_{\pizero} c}Q_{ce'}$, that is, the sum of row $c$ in $Q$, up to the diagonal; further, $\pi_2=a$, hence $s_2=\sum_{e'\myprec_{\pizero} a}Q_{ae'}=0$ (since there is no item before $a$ in $\pizero$), and so on. 
Below we show how to reconstruct $\pi$ from $\pizero$ and its code $s_{1:n-1}$.
\beq\label{eq:gmms-example}
\begin{array}{c|c|l|l|l}
	rank & s_i & \pi\text{ after insertion} & \pi_0\text{ after deletion} & \text{Factor in Eq}~\eqref{eq:gmms-like}\\\hline
	1& 2 &  [c] & [a,b,d,e,f] & \frac{\theta_1^2}{Z_1(\theta_1)}\\
	2& 0 &  [c,a] & [b,d,e,f] & \frac{\theta_2^0}{Z_2(\theta_2)}\\
	3& 1 &  [c,a,d] & [b,e,f] & \frac{\theta_3^1}{Z_2(\theta_3)}\\
	4& 1 &  [c,a,d,b] & [e,f] & \frac{\theta_4^1}{Z_3(\theta_4)}\\
	5& 0 &  [c,a,d,b,f] & [e] & \frac{\theta_5^0}{Z_4(\theta_5)}\\
	6& - &  [c,a,d,b,f,e]  & - & 1\\
\end{array}
\eeq
Finally, the probability of the ranking $\pi$ is the product of the factors in the last column, matching equation~\eqref{eq:gmms-like}.
\end{example}

\subsection{Entropy and KL-divergence between distributions}
\label{sec:bg-entropy}

Assume that the sample space is $\Omega$,  finite. For any distribution $P,P':\Omega\rightarrow [0,1]$ the {\bf
entropy} is defined as
\beq
H(P) \;=\;-\sum_{x\in \Omega} P(x) \ln P(x) \;=\;\expe_P[-\ln P]
\eeq
with the usual convention $0\ln 0=0$; w.l.o.g. the logarithms are in base $e$.

The function $H(P)$ is non-negative and concave on the space of all
distributions over $\Omega$. The minimum $H=0$ is attained for
deterministic distributions and the maximum $\ln|\Omega|$ is
attained for the uniform distribution over $\Omega$. 

The entropy measures the uncertainty in a given distribution. Closely
related to the entropy is the {\bf Kullbach-Leibler divergence}
between two distributions:
\beq
D(P||P') \;=\; \sum_{x\in\Omega} P(x) \ln \frac{P(x)}{P'(x)}
\eeq
The KL-divergence is asymmetric in $P,P'$. It is non-negative, convex
jointly in $(P, P')$ and attains the minimum $D(P||P') = 0$ iff $P\equiv P'$.

Note that the KL-divergence can be written as
\beq \label{eq:xent-general}
D(P||P') \;=\; \sum_{x\in\Omega} P(x) (\ln P(x)-\ln P'(x))
\;=\;\expe_P[-\ln P'] - H(P)
\eeq
The first term in the above expression is the {\em crossentropy} between $P$ and $P'$,
\beq
\xent(P||P')\;=\; \expe_P[-\ln P'].
\eeq
It is easy to see that $\xent(P||P')$ represents the log-likelihood of $P'$, viewed as a model, when the sampling distribution is $P$. The minimum of the crossentropy is attained when $P=P'$, the Maximum Likelihood model.
Moreover, we have $H(P)=\xent(P||P)$.

An {\em exponential family} is a family of probability distributions given by
\beq
P(x)\;=\;a(x)e^{\thetavec^Tt(x)-\ln Z(\thetavec)}
\eeq
where $a(x)$ and $t(x)$ are fixed functions, $\thetavec$ is a vector
of parameters, and $Z(\theta)=\expe_P[a(x)e^{\thetavec^Tt(x)}]$ is the
normalization constant. The functions $t(x)$ are called sufficient
statistics.

It is easy to see that the \gmms~model class is an exponential
familly, w.r.t. the parameters $\theta_{1:n-1}$.
\beq
P^{\gmms}(\pi)\;=\;e^{-\sum_{r=1}^{n-1}\ln \theta_r s_r(\pi|\sigma) - \sum_{r=1}^{n-1}\ln Z_{n-r+1}(\theta_r)}
\eeq
with $s_i(\pi|\sigma)$ being the sufficient statistics. More surprisingly, in \cite{MPhadnisPattersonBilmes:uai07} it has been shown that the (elements of) matrix $Q(\pi)$ are a sufficient statistics for any modal permutation $\sigma$, meaning that the entire set $\{P^{\gmms}_{\sigma,\thetavec}:\ssss_n\rightarrow (0,1]\}$ is an exponential family with parameters $\sigma$ and $\theta_{1:n-1}$.

For a standard exponential model family
  \beq
  p_\theta(X)\;=\;e^{\theta^TX-\psi(\theta)}
  \quad
  \text{with }X\in {\mathcal X} \text{ and } \theta\in \Theta,
  \eeq
  the entropy has the special form
  \beq
H(\theta)=\psi(\theta)-\theta^T \expe_\theta[X],
  \eeq
   the crossentropy is
  \beq
  \xent(\theta||\theta')\;=\;\psi(\theta')-(\theta')^T\expe_\theta[X]
  \eeq
  and the KL-divergence is
  \beq
  KL(\theta||\theta')\;=\;(\theta-\theta')^T\expe_\theta[X]-\psi(\theta)+ \psi(\theta').
  \eeq
  
\mmp{in intro -- la ce folosesc toate astea}

\section{Tractable computation of crossentropy in the \gmms~family}
\label{sec:crossent}
This section will introduce step by step a polynomial time algorithm for the computing the crossentropy $\xent(P\gmmsup_{\sigma,\thetavec}||P\gmmsup_{\sigma',\thetavec'})$. The existence of a tractable algorithm is a consequence of the close relationship between the \gmms~model and the inversion matrix $Q$. At first sight, computing the expectation in \eqref{eq:xent-general} for the \gmms~requires summation over a sample space of size $n!$. We bypass this challenge by showing that the crossentropy depends on matrix - vector products with submatrices of $Q$ which can be averaged tractably.

\subsection{Crossentropy as a function of sufficient statistics}
\label{sec:crossent-sbar}

We consider two distributions $P\gmmsup_{\idperm,\alpha_{1:n-1}},
P\gmmsup_{\sigma,\thetavec}$, respectively the sampling distribution
and the model, with $P,P'\in \gmms$.  We note that from equation
\eqref{eq:d-pisigma} and the definition of $Q$ it follows that w.l.o.g,
one of the modal permutations can be set to $\idperm$. We choose for
simplicity the set the modal permutation of the sampling distribution
to $\idperm$.
The cross-entropy between these can be written as follows.
\begin{prop}[Crossentropy as a function of average inversion counts]
\beq
\xent(P_{\idperm,\alphavec}||P_{\sigma,\thetavec})\;=\;
  \sum_{r=1}^{n-1}\left[\sbar_r(\sigma,\alpha_{1:r}) \ln \theta_r  + \ln Z_{n-r+1}(\theta_r)\right]
  \eeq
  where $\sbar_r=\expe_{\idperm,\alphavec}[s_r(\pi|\sigma)]$. 
\end{prop}

{\bf Proof}
 \beqa
\lefteqn{\xent(P_{\idperm,\alphavec}||P_{\sigma,\thetavec})}\\
  &=& E_{P_{\idperm,\alphavec}}\left[ -\ln P_{\sigma,\thetavec} \right]\\
  &=& \sum_{\pi\in\ssss_\E} P_{\idperm,\alphavec}(\pi) (-\ln P_{\sigma,\thetavec}(\pi))\\
&=& \sum_{\pi\in\ssss_\E} P_{\idperm,\alphavec}(\pi) \sum_{r=1}^{n-1}\left[  (s_r(\pi|\sigma) \ln \theta_r +\ln Z_{n-r+1}(\theta_r))\right]\\
&=&  \sum_{\pi\in\ssss_\E} P_{\idperm,\alphavec}(\pi) \sum_{r=1}^{n-1}\ln Z_{n-r+1}(\theta_r)+
\sum_{\pi\in\ssss_\E} P_{\idperm,\alphavec}(\pi) \sum_{r=1}^{n-1} (s_r(\pi|\sigma) \ln \theta_r\\
&=&   \sum_{r=1}^{n-1}\ln Z_{n-r+1}(\theta_r)+
\sum_{r=1}^{n-1} \ln \theta_r \left[\sum_{\pi\in\ssss_\E} P_{\idperm,\alphavec}(\pi)  s_r(\pi|\sigma)\right]\\
&=&   \sum_{r=1}^{n-1}\left[\sbar_r(\sigma,\alpha_{1:r}) \ln \theta_r  + \ln Z_{n-r+1}(\theta_r)\right]
\eeqa
\hfill $\Box$

This result is not surprising, being simply the expression of the (negative) log-likelihood of an exponential family model, separating the sufficient statistics from the continuous parameters. Notice that the discrete parameter $\sigma$ remains entangled with the sufficient statistc $\sbar_r$.  

We now focus on a single $\sbar_r$, representing the expected value of
$s_r(\pi|\sigma)$ under the sampling distribution
$\idperm,\alphavec$. Note that this term does not depend on the
parameters $\theta_{1:n}$. We will show  in Proposition \ref{prop:sbar-cond} that it is also independent of
$\alpha_{r+1:n-1}$.

Denote by $\E_r(\pi)=\{\pi_r,\ldots \pi_n\}$ the unordered set of
items left to be placed after ranks $1:r-1$ of $\pi$ have been
sampled. The codes $s_r(\pi|\idperm)$ and $s_r(\pi|\sigma)$ represent
respectively the position of $\pi_r$ in $\E_r(\pi)$, according to
$\idperm$, respectively $\sigma$. Hence,
$s_r(\pi|\sigma)=(\sigma\left|_{\E_r(\pi)})^{-1}(\pi_r)\right.$,
i.e. the rank of $\pi_r$ in $\E_r$ ordered by $\sigma$.  Hence,
$s_r(\pi|\idperm)=(\idperm\left|_{\E_r(\pi)})^{-1}(\pi_r)\right.$, i.e. the rank
of $\pi_r$ in $\E_r$  ordered by $\idperm$. For simplicity, let $k=n-r+1$. 

\begin{example}
Let $\pi=[c, a, d, b, f, e]$, $\idperm=[a,b,c,d,e,f]$ as before, and let $\sigma =[f,b,a,c,e,d]$. 
\beq\label{eq:srid-srsigma}
\begin{array}{c|c|c|c|l|l}
  \text{rank} & \pi_r & s_r(\pi_r|\idperm) & s_r(\pi_r|\sigma) & \idperm\text{ after deletion }& \sigma\text{ after deletion}\\
 \hline
	1 & c& 2 & 3 & [a,b,d,e,f] & [f,b,a,e,d]\\
	2 & a& 0 & 2 & [b,d,e,f]   & [f,b,e,d]\\
	3 & d& 1 & 4 & [b,e,f]     & [f,b,e]\\
	4 & b& 0 & 1 & [e,f]       & [f,e]\\
	5 & f& 1 & 0 & [e]         & [e]\\
	6 & e& - & - &  - & -\\
\end{array}
\eeq
\end{example}

\begin{prop}[Conditional expression for $\sbar_r$]
  \label{prop:sbar-cond}
For any $\sigma,\pi\in\ssss_n$, the expectation of $s_r(\pi|\sigma)$ under $P\gmmsup_{\idperm,\alphavec}$ is 
    \beq
  \sbar_r
\;=\;\frac{1}{Z_k(\alpha_r)}\sum_{\pi_{1:r-1}}P_{\idperm,\alpha_{1:r-1}}(\pi_{1:r-1})\left[
    \sum_{\pi_r\in \tilE_r} \alpha_r^{s_r(\pi|\idperm)} s_r(\pi|\sigma)\right]
\eeq
  \end{prop}
{\bf Proof} 
  \beqa
  \sbar_r
&=&\sum_{\pi\in\ssss_\E}\prod_{j=1}^{n-1}\frac{\alpha_j^{s_j(\pi|\idperm)}}{Z_{n-j+1}(\alpha_j)}
  s_r(\pi|\sigma)\\
  &=&\sum_{\pi_{1:r-1}}P_{\idperm,\alpha_{1:r-1}}(\pi_{1:r-1})\left[
  \sum_{\pi_r\in \tilE_r}\sum_{\ssss_{\tilE_{r+1}}}\frac{\alpha_r^{s_r(\pi|\idperm)}}{Z_k(\alpha_r)} \prod_{j=r+1}^{n-1}\frac{\alpha_j^{s_j(\pi|\idperm)}}{Z_{n-j+1}(\alpha_j)}
s_r(\pi|\sigma)\right]\\
  &=&\sum_{\pi_{1:r-1}}P_{\idperm,\alpha_{1:r-1}}(\pi_{1:r-1})\left[
    \sum_{\pi_r\in \tilE_r}\frac{\alpha_r^{s_r(\pi|\idperm)}}{Z_k(\alpha_r)}  s_r(\pi|\sigma)
\underbrace{\sum_{\ssss_{\tilE_{r+1}}}  \prod_{j=r+1}^{n-1}\frac{\alpha_j^{s_j(\pi|\idperm)}}{Z_{n-j+1}(\alpha_j)}}_{=1}\right]\\
  &=&\sum_{\pi_{1:r-1}}P_{\idperm,\alpha_{1:r-1}}(\pi_{1:r-1})\left[
    \sum_{\pi_r\in \tilE_r}\frac{\alpha_r^{s_r(\pi|\idperm)}}{Z_k(\alpha_r)}  s_r(\pi|\sigma)\right]
  \\
  &=&\frac{1}{Z_k(\alpha_r)}\sum_{\pi_{1:r-1}}P_{\idperm,\alpha_{1:r-1}}(\pi_{1:r-1})\left[
    \sum_{\pi_r\in \tilE_r} \alpha_r^{s_r(\pi|\idperm)} s_r(\pi|\sigma)\right]
  \eeqa
\hfill $\Box$

\begin{prop}[Matrix expression for $\sbar_r$] \label{prop:Qrup}
Let $A$ be a subset of $1:n$, and let $Q\uppper$ denote the lower triangular matrix obtained from $Q(\sigma)$ by zeroing out all the elements above the diagonal; $Q_{A,A}\uppper$ is the block $Q\uppper_{i,i'}$ with $i,i'\in A$, $P(A)=Pr[ \pi_{1:r-1}\in \E\setminus A\,|\,\idperm, \alpha_{1:r-1}]$ and $\bbone_k$ is the column vector of all ones, of dimension $k$, and $k=n-r+1$ as before. Then,
  \beq\label{eq:sbar-QP} 
  \sbar_r\;=\;\frac{1}{Z_{n-r+1}(\alpha_r)}\bar{\alpha}_r^T\sum_{A\subset \E, |A|=k}
\!\!\!\!\!
  Q\uppper_{A,A}P(A)\bbone_k
  \;=\;
\frac{1}{Z_{n-r+1}(\alpha_r)}\bar{\alpha}_r^T\qbar_k\uppper(\sigma,\alpha_{1:r-1})\bbone_k
  \eeq
  with
  \beq
  \bar{\alpha}_r^T\;=\;[\,1\,\alpha_r\,\alpha_r^2\,\ldots\,\alpha_r^{k-1}\,],
 \eeq
\end{prop}

{\bf Proof} At step $r$, $s_r$ is sampled from the set $A = \{1:n \}
\setminus \{\pi_{1:r-1}\}=\{a_1,\ldots\,a_k\}$. Each column $i$ of $Q\uppper_{A,A}$ contains
1's for the elements $a_1,\ldots a_{i-1}$ that are before $a_i$ in $\sigma$. Thus, if $a_i$ is
sampled, $s_r(i)=\sum_{j=1}^{n} (Q_{A,A}\uppper)_{ij}$. \mmp{recheck this} The vector of
these values is obtained by $Q_{A,A}\uppper\bbone_k$. Under the
sampling distribution, the item at rank $i$ has probability
$\alpha_r^{i-1}/Z_{n-r+1}(\alpha_r)$ to be sampled. This probability
is the same for any set $A$, since it depends only on the sampled
rank,not on the item at that rank. Finally, $P(A)$ is the probability
that after $r-1$ steps we are left with set $A$. \hfill $\Box$
 
Hence, $\qbar_k\uppper$ is the average under the sampling distribution of all principal blocks of size $k$ in $Q\uppper$. 

Up to here, only the model needed to be in the \gmms~class, while the sampling distribution can be any stagewise sampling model\footnote{Called a {\em full} model by \cite{fligner:86}.}

\subsection{Entropy and KL divergence}
\label{sec:xent-entropy}
If $\sigma=\idperm$, then $Q(\idperm)$ is a upper triangular, with $1$ strictly above the diagonal. Thus $Q\uppper=Q$ and any $Q\uppper_{A,A}$ will have the same pattern. 
  Therefore
  $(\qbarl_k)_{ij}=1$ if $j>i$ and 0 otherwise, and
\beq
\sbar_r\;=\;\frac{\sum_{i=1}^{n-r}i\alpha_r^i}{Z_{n-r+1}(\alpha_r)}
\eeq
  Hence, the entropy of a \gmms~ is independent of the location parameter $\sigma$ and equals
  \beq
  H(\thetavec)\;=\;
  \sum_{r=1}^{n-1}\left[\frac{\sum_{i=1}^{n-r}i\theta_r^i}{Z_{n-r+1}(\theta_r)}  \ln \theta_r  + \ln Z_{n-r+1}(\theta_r)\right]
  \eeq
  This can be expressed also as
  \beq
  H(\thetavec)\;=\;
  \sum_{r=1}^{n-1}\left[\expe_{\geom(\theta_r,n-r+1)}[i]  \ln \theta_r  + \ln Z_{n-r+1}(\theta_r)\right]
  \eeq
  which is the expression of the entropy in an exponential family, namely the entroypy of $n-1$ independent geometric distributions.

  The KL-divergence is the difference between cross-entropy and entropy, hence

  \begin{prop}[KL-divergence for \gmms]
 The KL-divergence between two \gmms~models with parameters
 $\sigma,\thetavec$ and $\sigma',\thetavec'$ is
 \beq \label{eq:kl}
 KL(\,P_{\sigma,\thetavec}|| P_{\sigma',\thetavec'})\;=\;\sum_{r=1}^{n-1}\frac{1}{Z_{n-r+1}(\theta_r)}\left[\bar{\theta}_r\qbarl_{n-r+1}(\sigma(\sigma')^{-1},\theta_{1:r})\bbone_{n-r+1}\ln \frac{\theta'_r}{\theta_r} + \ln \frac{Z_{n-r+1}(\theta'_r)}{Z_{n-r+1}(\theta_r)}\right].
 \eeq
\end{prop}

\section{Calculation of $\qbarl_k(\sigma, \alpha_{1:r-1})$}
\label{sec:qbar}
It would appear from Proposition \ref{prop:Qrup} that to obtain $\qbar_r\uppper$, one
would have to compute a sum with $\binom{n}{k}$ terms, for $k=2,\ldots
n-1$, that is, to sum $\bigOO(2^n)$ matrices. However, by considering
the recursive dependencies between $\qbarl_{r-1}$ and $\qbarl_r$, we
can obtain all $\qbar\uppper_{1:n-1}$ in polynomial time and memory,
as described in algorithm \ref{alg:qbar}.

\begin{prop}
Algorithm \ref{alg:qbar} calculates $\qbar\uppper_{1:n-1}$ in polynomial time.
\end{prop}

{\bf Proof idea} It is useful in this section to assume that the set
$\E=1:n$, hence its items are numbers in $1:n$, and $\idperm$ is the
natural order. The idea is to note that, for any rank $r$, the
contributions to $\qbarl_r$ can come only from pairs $a,b\in \E$,
$a<b$, for which $Q_{ba}=1$. In other words, from the inversions of
$\sigma$ w.r.t. $\idperm$. Hence, instead of enumerating sets $A$, we
enumerate inverted pairs $a,b$, and track where they can appear in
$\qbarl_r$. Since there are at most $n(n-1)/2$ inversions, and $O(n)$
steps and ranks, the algorithm is polynomial.

{\bf Proof} 
Denote by $\pijtab^r_{ab,ij}$ the following probability, with $a,b\in \E,\,a<b,1\leq i<j\leq k$.
\beq\label{eq:pabii'}
\pijtab^r_{ab,ij}\;=\;P[\text{after $r-1$ samples removed from }1:n,\, \rank(a)=i, \,\rank(b)=j],
\eeq
or, in other words, denoting by $E_{r-1,ab,ij}$ the event
$E_{r-1,ab,ij}\;=\;\{\,(\idperm\left|\right._{[n]\setminus   \pi_{1:r-1}})=a,\, (\idperm\left|\right._{[n]\setminus   \pi_{1:r-1}})_{j}=b\,\}$, we denote $\pijtab^r_{ab,ij}=Pr[ E_{r-1,ab,ij'}]$. For each rank $r$, we also maintain a list $\listab^r_{ab}$ of all the ranks  $(i,j)$ where the pair $(a,b)$ can appear in $\qbarl_r$. Since $b>a$, we only need to consider locations  with $j>i$. All probabilities in this proof are under $P\gmmsup_{\idperm,\alphavec}$, with $\sigma$ a fixed permutation.

\paragraph{Initialization} At rank $r=1$, $k=n$, and $(a,b)$ can occupy only 1 position, $(i,j) = (a,b)$. One can verify that $\sbar^1=\sum_{i=1}^n\alpha_1^{i-1}(\sigma^{-1}(i)-1)
\:=\; \sum_{i=1}^n\alpha_1^{i-1}\sum_{j}Q_{ij}\uppper\;=\;\bar{\alpha}_1Q\uppper\bbone_n$

\paragraph{Recursion} 
Assuming  $\pijtab_{ab}^{r-1},\listab_{ab}^{r-1}$ for the
previous rank are given, we have to enumerate all the positions $(i_0,j_0)$ and
their respective probabilities, that can lead to the event
$E_{ab,ij}$. If $\rank a =i_0$, and $\rank b=j_0$ before the
$(r-1)$-th sample, then after sampling an item, these ranks can only
stay the same decrease. Namely, $i\gets i_0,j\gets j_0$ if the
deleted element is after $b$, $i\gets i_0,j\gets j_0-1$ if the
deleted element is between $a$ and $b$, and $i\gets i_0-1,j\gets
j_0-1$ if the selected element is before $a$. The events are encoded
in the values of the code $s_{r}$ and their probabilities are given by
the $\operatorname{geom}(\alpha_r, n-r+1)$ distribution. Therefore,
from $\listab_{ab}^{r-1}$ we can get all $(i,j)$ positions of $a,b$ in
$\listab_{ab}^r$ and their probabilities. Note that there are multiple
ways for $(a,b)$ to arrive in position $(i,j)$ if
$\listab_{ab}^{r-1}$ has more than one entry. Therefore, each time a
position $(i,j)$ is enumerated, its probabilty is added to the current
value for $\pijtab_{ab,ij}^r$. This addition corresponds to the total
probability of $E_{ab,ij}$ as a sum over mutually exclusive events.

Finally, in steps \ref{step:add-ab-to-Q}--\ref{step:end-add-ab-to-Q}, the contributions of item pair $(a,b)$ are added to $\qbarl_r$. 
\hfill $\Box$

\begin{algorithm}[h]
  \caption{\label{alg:qbar}Recursive computation of $\qbarl_{1:n-1}$}
  \begin{algorithmic}[1]
    \State {\bf Input} Inversion matrix $Q$, parameters $\alpha_{1:n-1}$
    \State {\bf Initialize}
    \State $\qbarl_1\gets Q\uppper$
    \State $\listab_{ab}^1\gets [(a,b)]$, $\pijtab^1_{ab,ab}\gets Q_{ab}\uppper,\,$
    for $1\leq a < b\leq n$ and $Q_{ba}=1$. \label{step:init-pij-list}
     
    \For{$r=1:n-1$}
        \State $k\gets n-r$ 
        \State $\qbarl_r\gets {\mathbf 0}_{k\times k}$
        \For{ $a,b$ with $1\leq a < b\leq n$  and $Q_{ba}=1$}
            \State $\listab_{ab}^r\gets []$
            \State $\pijtab_{ab}^r \gets []$  
            \For{$(i_0,j_0)\in \listab_{ab}^{r-1}$}
                \If{$j_0\leq k$} \Comment{rank $r$ sampled after $j_0$, $s_r\geq j_0$}
                   \State $i\gets i_0$, $j\gets j_0$
                   \State $\operatorname{insert}(\listab^r_{ab}, (i,j))$ \Comment{if not already existing}
                   \State $\pijtab^{r}_{ab,ij} \gets \pijtab^r_{ab,ij}+  \pijtab^{r-1}_{ab,i_0j_0}\sum_{s=j}^{k-1}\alpha_r^s/Z_{k}(\alpha_r)$
                 \EndIf
                \If{$i_0+1<j_0$} \Comment{rank $r$ sampled between $i_0$ and $j_0$, $i_0\leq s_r< j_0$}
                   \State $i\gets i_0$, $j\gets j_0-1$
                   \State $\operatorname{insert}(\listab^r_{ab}, (i,j))$
                   \State $\pijtab_{ab,ij}^r \gets \pijtab_{ab,ij}^r
                  +\pijtab_{ab,i_0j_0}^{r-1}\sum_{s=i_0}^{j_0-1}\alpha_r^s/Z_{k}(\alpha_r)$
                 \EndIf
                \If{$i_0 >1$} \Comment{rank $r$ sampled before $i_0$, $s_r< i_0$}
                   \State $i\gets i_0-1$, $j\gets j_0-1$
                   \State $\operatorname{insert}(\listab^r_{ab}, (i,j))$
                   \State $\pijtab_{ab,ij}^r \gets \pijtab_{ab,ij}^r+\pijtab_{ab,i_0j_0}^{r-1}\sum_{s=0}^{i_0-1}\alpha_r^s/Z_{k}(\alpha_r)$
                 \EndIf
                \For{$l = 1:\operatorname{length}(\listab^r_{ab})$} \Comment{add $(a,b)$ contributions to $\qbarl$} \label{step:add-ab-to-Q}
                \State $(i,j)\gets \listab^r_{ab}(l)$
                \State $\qbarl_{r,ji}\gets \qbarl_{r,ji}+\pijtab_{ab,ij}^r$ \label{step:add-ab-to-Q}
                \EndFor \label{step:end-add-ab-to-Q}
            \EndFor
        \EndFor
    \EndFor
  \State {\bf Output} $\qbar^1,\,\ldots \qbar^{n-1}$
  \end{algorithmic}
  \end{algorithm}

\begin{corollary}
  $\qbar_r=\sum_{A\subset \E, |A|=k}Q_{A,A}P(A)$, with $k=n-r+1$, can be calculated by a modification of Algorithm \ref{alg:qbar} that is initialized with $Q$ and tracks all item pairs, both above and below the diagonal. Specifically
  \bit
\item $\qbarl_r\gets \qbar_r$ everywhere in Algorithm \ref{alg:qbar}.
\item In  step \ref{step:init-pij-list} initializations are done for all pairs $(a,b)$ regardless of the value of $Q_{ab}$.
\item Step \ref{step:add-ab-to-Q} becomes
  \[ \left\{\begin{array}{ll} \qbar_{r,ij}\gets \qbar_{r,ij}+\pijtab^r_{ab,ij}
  & \text{\bf if } Q_{ab}=1\\
  \qbar_{r,ji}\gets \qbar_{r,ji}+\pijtab^r_{ab,ij} 
  & \text{\bf if } Q_{ba}=1\\
  \end{array} \right.
  \]
  \eit
\end{corollary}

{\bf Proof idea} For every pair $a<b$ in $\E$ exactly one of $Q_{ab},Q_{ba}$ equals 1. Thus, the data structures $\listab_{ab}^r,\pijtab_{ab}^r$ are associated to the $Q$ entry that is non-zero. Otherwise, the recursion proceeds identically until contribution of pair $(a,b)$ is added to $\qbar_r$, when again we must check which of $Q_{ab},Q_{ba}$ is non-zero. This is possible because the zero elements of $Q$ no not contribute to any element of $\qbar_r$. \hfill $\Box$

\section{Examples and experiments}
\label{sec:exe}

\begin{example}[KL-divergence for same central permutation]
  Let $\sigma=\sigma'=\idperm$ (w.l.o.g), while the parameters $\thetavec\neq\thetavec'$. Then, the KL-divergence simplifies to
\beq
 KL(\,P_{\sigma,\thetavec}|| P_{\sigma,\thetavec'})\;=\;  \sum_{r=1}^{n-1}\left[\expe_{\geom(\theta_r,n-r+1)}[i]  \ln\frac{\theta'_r}{\theta_r}  + \ln \frac{Z_{n-r+1}(\theta'_r)}{Z_{n-r+1}(\theta_r)}\right],
  \eeq
 representing the sum of KL-divergences between geometric distributions on $\{0,\ldots n-r+1\}$.   Moreover, the terms for which $\theta_r=\theta'_r$ are zero, hence we can write the above sum as
\beq
 KL(\,P_{\sigma,\thetavec}|| P_{\sigma,\thetavec'})\;=\;  \sum_{r:\theta_r\neq\theta_r'}\left[\expe_{\geom(\theta_r,n-r+1)}[i]  \ln\frac{\theta'_r}{\theta_r}  + \ln \frac{Z_{n-r+1}(\theta'_r)}{Z_{n-r+1}(\theta_r)}\right].
  \eeq

\end{example}

\begin{example}[KL-divergence for same parameters, different central permutations]
  Let $\sigma\neq \sigma'$, while the parameters $\thetavec=\thetavec'$. Then, the KL-divergence simplifies to
\beq
 KL(\,P_{\sigma,\thetavec}|| P_{\sigma',\thetavec})\;=\;  \sum_{r=1}^{n-1}(\sbar_r(\sigma(\sigma')^{-1},\theta_{1:r})-\expe_{\geom(\theta_r,n-r+1)}[i])  \ln \theta_r
  \eeq
From the above, we infer immediately that $\sbar_r(\sigma,\theta_{1:r})\geq \sbar_r(\idperm,\theta_{1:r})$ for all $\theta_{1:r}\in [0,1]$ and all $\sigma\in\ssss_n$.
 
\end{example}

\section{Discussion}
\label{sec:discussion}
Once again, the \gmms~model has proved its versatiliy. This paper showed
how to separate integration over the continuous spread parameters
$\theta$ from summation over permutations in $\ssss_n$, and how the
existence of the $Q$ representation alllows for the later to be
performed in polynomial time. This work opens the way for advancing
the study of exponential families over discrete parameter spaces, in
order to expand to them the advantages and the intuitition that
standard exponential families enjoy \cite{Zhang2014}.

We note that for the \gmms$^v$ model, as well as for the Recursive Inversion Models \cite{WagnerM:RIM22}, which are essentially extensions thereof, we currently do not know if tractable algorithms exist for the computation of the entropy and crossentropy.

\section*{Acknowledgements}
The author thanks the organizers of the Interactions of Statistics and
Geometry (ISAG) II Workshop, and especially to Jun Zhang, for the interesting
presentations that are the impetus for this work. The author
acknowledges partial support from the NSF MMS 2019901 award as well as
the support of the Institute for Mathematical Sciences at the National
University of Singapore, and of the Datashape group at INRIA Saclay. 

\bibliographystyle{apalike}
\bibliography{refs}
\end{document}